\documentclass[a4paper,12pt]{article}
\usepackage{cmap}                        
\usepackage[cp1251]{inputenc}            
\usepackage[english]{babel}
\usepackage[left=2cm,right=2cm,top=2cm,bottom=2cm]{geometry} 
\usepackage{amssymb}
\usepackage{amsmath}

\begin{document}

\begin{center}\Large
\textbf{ON GENERALIZATIONS OF BAER'S THEOREMS  
ABOUT THE HYPERCENTER OF A FINITE GROUP}
\end{center}
\begin{center} V. I. Murashka \end{center}
\begin{center} \{mvimath@yandex.net\}\end{center}
\begin{center} Francisk Skorina Gomel State University, Gomel\end{center}

\medskip

\textbf{Abstract.} We investigate the intersection of normalizers and $\mathfrak{F}$-subnormalizers of different types of systems of subgroups ($\mathfrak{F}$-maximal, Sylow, cyclic primary).  We described all formations $\mathfrak{F}=\underset{i\in I}\times\mathfrak{F}_{\pi_i}$ for which the intersection of normalizers of all $\mathfrak{F}_i$-maximal subgroups of   $G$ is the $\mathfrak{F}$-hypercenter of $G$ for every group $G$.  Also we described all formations $\mathfrak{F}$   for which the intersection of $\mathfrak{F}$-subnormalizers of all Sylow (cyclic primary) subgroups of   $G$ is the $\mathfrak{F}$-hypercenter of $G$ for every group $G$.

\textbf{Keywords:} saturated formation, hereditary formation, $\mathfrak{F}$-hypercenter, $\mathfrak{F}$-subnormalizer, intersection of subgroups.

\textbf{Mathematic Subject Classification}(2010): 20D25,  20F17, 20F19.

\section{Introduction}

All considered groups are finite. In \cite{h1} R. Baer showed that from one hand the hypercenter $\mathrm{Z}_\infty(G)$ of a  group $G$ coincides with the intersection of all maximal nilpotent subgroups of $G$ and from another hand $\mathrm{Z}_\infty(G)$ coincides with the intersection of normalizers of all Sylow subgroups  of $G$.

The concept of hypercenter was extended on classes of groups (see \cite[p. 127--128]{s6} or \cite[p. 6--8]{s5}). Let $\mathfrak{X}$ be a class of groups. A chief factor $H/K$ of a group $G$ is called   $\mathfrak{X}$-central if $(H/K)\leftthreetimes G/C_G(H/K)\in\mathfrak{X}$. A normal subgroup $N$ of $G$ is said to be $\mathfrak{X}$-hypercentral in $G$ if $N=1$ or $N\neq 1$ and every chief factor of $G$ below $N$ is $\mathfrak{X}$-central. The $\mathfrak{X}$-hypercenter $\mathrm{Z}_\mathfrak{X}(G)$ is the product of all normal   $\mathfrak{X}$-hypercentral subgroups of $G$. So if $\mathfrak{X}=\mathfrak{N}$ is the class of all nilpotent groups then $\mathrm{Z}_\infty(G)=\mathrm{Z}_\mathfrak{N}(G)$ for every group $G$.

In \cite{h2} A.\,V. Sidorov  showed that for a soluble group $G$ the intersection of all maximal subgroups of nilpotent length at most $r$ is $\mathrm{Z}_{\mathfrak{N}^r}(G)$.   Beidleman  and Heineken \cite{h3} studied the properties of the intersection $\mathrm{Int}_\mathfrak{F}(G)$ of $\mathfrak{F}$-maximal subgroups of a group $G$ in case when $G$ is soluble and $\mathfrak{F}$ is a hereditary saturated formation.

Let $F$ be the canonical local definition of a local formation $\mathfrak{F}$. Then $\mathfrak{F}$ is said to satisfy the boundary condition \cite{h4} if  $\mathfrak{F}$ contains every group $G$  whose all maximal subgroups belong to $F(p)$ for some prime $p$.

 A.\,N. Skiba \cite{h4} showed that the equality $\mathrm{Int}_\mathfrak{F}(G)=\mathrm{Z}_\mathfrak{F}(G)$ holds for every group $G$ if and only if  a hereditary saturated formations $\mathfrak{F}$ satisfies the boundary condition. This and further results was included in the first chapter of \cite{s5}.

The intersection of normalizers of different systems of subgroups is the main theme of many papers. In \cite{h5} Baer  considered the intersection of normalizers of all subgroups of a group.  Wielandt \cite{h6} studied  the intersection of normalizers of all subnormal subgroups of a group. Li and Shen \cite{h7} considered the  intersection of normalizers of all derived subgroups of all subgroups of a  group. 

Let $\sigma=\{\pi_i | i\in I\}$ be a partition of $\mathbb{P}$ into disjoint subsets, $\mathfrak{X}_i$ be a class of groups such that $\pi(\mathfrak{X}_i)=\pi_i$.   Then  $\underset{i\in I}\times\mathfrak{X}_{\pi_i}=(G=\underset{i\in I}\times \mathrm{O}_{\pi_i}(G)| \mathrm{O}_{\pi_i}(G)\in\mathfrak{X}_i)$.
Recall that $\mathfrak{G}_\pi$ is the class of all $\pi$-groups.  Hence $\mathfrak{N}=\underset{p\in \mathbb{P}}\times\mathfrak{G}_{p}$.

In \cite{h8} author showed that if $\mathfrak{F}=\underset{i\in I}\times\mathfrak{G}_{\pi_i}$ then for any group $G$ the intersection of all normalizers of all $\pi_i$-maximal subgroups of  $G$ for all $i\in I$ coincides with the $\mathfrak{F}$-hypercenter. So the general problem is

\textbf{Problem A.} Let $\Sigma$ be a subgroup functor. What can be said about the intersection of normalizers of subgroups from $\Sigma(G)$?

Recall \cite[p. 206]{s9} that a subgroup functor is a function $\tau$ which
 assigns to each group $G$ a possibly empty set $\tau(G)$ of subgroups of $G$ satisfying
 $f(\tau(G))=\tau(f(G))$ for any isomorphism $f: G\rightarrow G^*$.

 \textbf{Definition 1.} Let $\mathfrak{X}$ be a class of groups and $G$ be a group. Then $\mathrm{NI}_\mathfrak{X}(G)$ is the intersection of all normalizers of $\mathfrak{X}$-maximal subgroups of $G$.

  The following proposition shows that if $\mathfrak{F}$ is a hereditary saturated formation and $\pi(\mathfrak{F})=\mathbb{P}$  then the equality $\mathrm{NI}_\mathfrak{F}(G)=\mathrm{Int}_\mathfrak{F}(G)$ holds for every group $G$.

\textbf{Proposition 1.}  \emph{Let $\mathfrak{F}$ be a hereditary saturated formation and $\pi=\pi(\mathfrak{F})$. Then \linebreak $\mathrm{O}^{\pi'}(\mathrm{NI}_\mathfrak{F}(G))= \mathrm{Int}_\mathfrak{F}(G)$ for every group $G$. }

The following theorem generalizes two above mentioned   Baer's theorems about the hypercenter:

\textbf{Theorem A.} \emph{Let  $\sigma=\{\pi_i | i\in I\}$ be a partition of $\mathbb{P}$ into disjoint subsets and $\mathfrak{F}_i$ be a hereditary saturated formation such that $\pi(\mathfrak{F}_i)=\pi_i$ and $\mathfrak{F}=\underset{i\in I}\times\mathfrak{F}_{i}$. The following statements are equivalent}:

$(1)$ \emph{$\mathfrak{F}_i$  satisfies the boundary condition in the universe of all $\pi_i$-groups for all $i\in I$. }

$(2)$ \emph{For every group $G$ holds $\underset{i\in I}\bigcap\mathrm{NI}_{\mathfrak{F}_i}(G)=\mathrm{Z}_\mathfrak{F}(G)$. }

\textbf{Corollary A.1 \cite{h1}.} \emph{The hypercenter of a group $G$ is the intersection of all normalizers of all Sylow subgroups of $G$.}

\textbf{Corollary A.2 \cite{h8}.} \emph{Let  $\sigma=\{\pi_i | i\in I\}$ be a partition of $\mathbb{P}$ into disjoint subsets, $\mathfrak{F}=\underset{i\in I}\times\mathfrak{G}_{\pi_i}$ and $G$ be a group. Then   the intersection of all normalizers of all $\pi_i$-maximal subgroups of  $G$ for all $i\in I$ is the $\mathfrak{F}$-hypercenter of $G$.}

From proposition 1 and theorem A when $|I|=1$ it follows that our theorem A   extends   theorem A from \cite{h4}:

\textbf{Corollary A.3.} \emph{Let $\mathfrak{F}$ be a hereditary saturated formation and $\pi(\mathfrak{F})=\mathbb{P}$. The equality $\mathrm{NI}_\mathfrak{F}(G)=\mathrm{Z}_\mathfrak{F}(G)=\mathrm{Int}_\mathfrak{F}(G)$ holds for every group $G$ if and only if  $\mathfrak{F}$ satisfies the boundary condition. }

\textbf{Corollary A.4 \cite{h1}.} \emph{The hypercenter of a group $G$ is the intersection of all maximal nilpotent subgroups  of $G$.}

Let $\mathfrak{X}$ be a class of groups. Recall that a subgroup $H$ of a group $G$ is called $\mathfrak{F}$-subnormal if either $H=G$ or there is a maximal chain of subgroups $H=H_0<H_1<\dots< H_n=G$ such that $H_i/\mathrm{Core}_{H_i}(H_{i-1})\in\mathfrak{X}$ for all $i=1, \dots, n.$

Let $\mathfrak{X}$ be a class of groups. A $\mathfrak{X}$-subnormalizer \cite[p. 380]{s8} of a subgroup $H$ of a group $G$ is a subgroup $T$ of $G$ such that $H$ is $\mathfrak{X}$-subnormal in $T$ and if  $H$ is $\mathfrak{X}$-subnormal in $M$ and $T\leq M$ then $T=M$. It is clear that a $\mathfrak{X}$-subnormalizer always exists but may be not unique.

\textbf{Problem B.} Let $\Sigma(G)$ be a subgroup functor and $\mathfrak{F}$ be a formation. What can be said about the intersection $\mathrm{SI}_\Sigma^\mathfrak{F}(G)$ of $\mathfrak{F}$-subnormalizers of subgroups from $\Sigma(G)$?

If $\Sigma(G)$ is the set of all maximal subgroups of $G$ then this intersection coincides with $\Delta_\mathfrak{F}(G)$ where $\Delta_\mathfrak{F}(G)$ is the intersection of all $\mathfrak{F}$-abnormal maximal subgroups of $G$. According to \cite[p. 96]{s7} if $\mathfrak{F}$ is a hereditary saturated formation then $\Delta_\mathfrak{F}(G)/\Phi(G)=\mathrm{Z}_\mathfrak{F}(G/\Phi(G))$.

\textbf{Proposition 2.}\emph{ Let $\mathfrak{F}$ be a hereditary formation and $\Sigma$ be a subgroup functor. Then $\mathrm{SI}_\Sigma^\mathfrak{F}(G)$ is the product of normal subgroups $N$ of a group $G$ such that $H$ is $\mathfrak{F}$-subnormal in $HN$ for every $H\in\Sigma(G)$.   }

A.F. Vasil'ev and T.I. Vasil'eva \cite{h9} studied a class of groups $w\mathfrak{F}$  whose all Sylow subgroups are $\mathfrak{F}$-subnormal for a given hereditary saturated formation $\mathfrak{F}$. Let us note that in this case $\mathrm{Z}_{w\mathfrak{F}}(G)$ lies in the intersection of all $\mathfrak{F}$-subnormalizers of all Sylow subgroups of a group $G$. Author \cite{h10} studied a class of groups $v\mathfrak{F}$  whose all cyclic primary subgroups are $\mathfrak{F}$-subnormal for a given hereditary saturated formation $\mathfrak{F}$. Again $\mathrm{Z}_{v\mathfrak{F}}(G)$ lies in the intersection of all $\mathfrak{F}$-subnormalizers of all cyclic primary subgroup of a group $G$.

In this paper we count the unit group as cyclic primary subgroup and also as Sylow subgroup.

\textbf{Theorem B.} \emph{Let $\mathfrak{F}$ be a hereditary saturated formation. The following statements are equivalent:}

$(1)$ \emph{There exists a partition $\sigma=\{\pi_i | i\in I\}$ of $\pi(\mathfrak{F})$ into disjoint subsets such  that $\mathfrak{F}=\underset{i\in I}\times\mathfrak{G}_{\pi_i}$.}

$(2)$ \emph{The intersection of all $\mathfrak{F}$-subnormalizers of all cyclic primary subgroups of  $G$ is the $\mathfrak{F}$-hypercenter of $G$ for every group $G$}.

$(3)$ \emph{The intersection of all $\mathfrak{F}$-subnormalizers of all Sylow subgroups of  $G$ is the $\mathfrak{F}$-hypercenter of $G$ for every group $G$}.

Note that in the universe of all soluble groups the concepts of a subnormal subgroup and a $\mathfrak{N}$-subnormal subgroup coincides. It is well known that if a Sylow subgroup $P$ of $G$ is subnormal in $G$ then it is normal in $G$. Hence a $\mathfrak{N}$-subnormalizer of a Sylow subgroup $P$  of a soluble group $G$ is just the normalizer of $P$ in $G$.  So theorem B can be viewed as the generalization   of R. Baer's theorem about the intersection of normalizers of Sylow subgroups.

\emph{Remark.} Formations $\mathfrak{F}=\underset{i\in I}\times\mathfrak{G}_{\pi_i}$ are lattice formations, i.e. formations were $\mathfrak{F}$-subnormal subgroups form a sublattice of the subgroup's lattice of every group (for example see chapter 6.3 of \cite{s9}). Also properties of the $\mathfrak{F}$-hypercenter and the $\mathfrak{F}$-residual for such formations was studied by author in \cite{h8}. A.\,N. Skiba extends the theory of nilpotent groups on such classes (for example see \cite{h11}).

\section{Preliminaries}

We use standard notation and terminology that if necessary can be found in \cite{s8}. Recall some of them that are important in this paper. By  $\mathbb{P}$ is denoted the set of all primes; $\pi(G)$ is the set of all prime divisors of the order of $G$; $\pi(\mathfrak{F})=\underset{G\in\mathfrak{F}}\cup\pi(G);$  a group  $G$ is called $\pi$-group if $\pi(G)\subseteq\pi$;    $Z_p$ is the cyclic group of order $p$;      $\mathrm{O}_\pi(G)$ is the greatest normal $\pi$-subgroup $G$; $\mathrm{O}^\pi(G)$ is the smallest subgroup of $G$ such that $\pi(G/\mathrm{O}^\pi(G))\subseteq\pi$ ;   $G'$ is the derived subgroup of $G$;  $G^\mathfrak{F}$ is the $\mathfrak{F}$-residual for a formation $\mathfrak{F}$;      $\mathrm{O}_{p', p}(G)$ is the    $p$-nilpotent radical of $G$ for $p\in\mathbb{P}$, it also can be defined by $\mathrm{O}_{p', p}(G)/\mathrm{O}_{p'}(G)=\mathrm{O}_p(G/\mathrm{O}_{p'}(G))$; $\Phi(G)$ is the Frattini subgroup of a group $G$;  $G=N\leftthreetimes M$ is the semidirect product of groups $M$ and $N$  ($N\triangleleft G$ and $N\cap M=1$);
 $\mathfrak{G}_\pi$ ($\mathfrak{S}_\pi$, $\mathfrak{N}_\pi$) is the class of  (soluble, nilpotent) $\pi$-groups, where $\pi\subseteq\mathbb{P}$.

A class of groups $\mathfrak{F}$ is called a formation if from $G\in\mathfrak{F}$ and $N\triangleleft G$
 it follows that $G/N\in\mathfrak{F}$ and from  $H/A\in \mathfrak{F}$ and $H/B\in \mathfrak{F}$ it
 follows that $H/A\cap B\in \mathfrak{F}$.

 A class of groups $\mathfrak{X}$ is called hereditary  if from $G\in\mathfrak{X}$ and
 $H\leq G$ it follows that $H\in\mathfrak{X}$.

 Let $\mathfrak{F}$ and $\mathfrak{K}$ be formations then $\mathfrak{FK}=(G| G^\mathfrak{K}\in\mathfrak{F})$ is also formation.

A class of groups  $\mathfrak{X}$ is called saturated  if from $G/\Phi(G)\in\mathfrak{X}$ it
 follows that $G\in\mathfrak{X}$.

By well known Gashutz-Lubeseder-Shmid Theorem saturated formations are exactly local formations,
i.e. formations $\mathfrak{F}=LF(f)$ defined by a formation function $f$: $LF(f)=\{
G\in\mathfrak{G}|$ if $H/K$ is a chief factor of $G$ and $p\in\pi(H/K)$ then $G/C_G(H/K)\in
f(p)\}$. Among all possible local definitions of a local formation $\mathfrak{F}$ there is exactly
one, denoted by $F$, such that $F$ is integrated ($F(p)\subseteq \mathfrak{F}$ for all
$p\in\mathbb{P}$)   and full ($\mathfrak{N}_pF(p)=F(p)$ for all $p\in\mathbb{P}$). $F$ is called
the canonical local definition of $\mathfrak{F}$.

Let $\mathfrak{F}$ be a local formation, $F$ be its canonical local definition  and $G$ be a group. Then a chief factor $H/K$ of a group $G$ is $\mathfrak{F}$-central if and only if $G/C_G(H/K)\in F(p)$ for all $p\in\pi(H/K)$ (see \cite[p.  6]{s5}).

Let $\mathfrak{F}=LF(F)$ be a hereditary local formation, $F$ be its canonical local definition  and $\pi=\pi(\mathfrak{H})$. Then $\mathfrak{F}$ is said to  satisfy the boundary condition in the universe of all $\pi$-groups if $\mathfrak{F}$ contains every $\pi$-group  whose all maximal subgroups belong to $F(p)$ for some prime $p$.

The following lemma can be found in \cite[p. 239]{s7}. For reader's convenience, we give a direct proof.

\textbf{Lemma 2.1.} \emph{Let $\mathfrak{X}$ be a saturated homomorph and $N$ be a normal subgroup of a group $G$. Then for every $\mathfrak{X}$-subgroup $H/N$ of $G/N$ there exists a $\mathfrak{X}$-subgroup $M$ of $G$ such that $MN/N=H/N$.}

\emph{Proof.} Let $H/N$ be a $\mathfrak{X}$-subgroup of $G/N$. Let us show that there exists a $\mathfrak{X}$-subgroup $K$ of $G$ such that $KN/N=H/N$. Let $M$ be a minimal subgroup of $H$ such that $MN=H$ (i.e if $M_1<M$ then $M_1N<H$). Assume that there is a maximal subgroup $M_1$ of $M$ such that $M_1(M\cap N)=M$. Then $M_1N=H$, a contradiction. Hence $M\cap N\leq\Phi(M)$. Since $\mathfrak{X}$  is saturated and $H/N=MN/N\simeq M/M\cap N\in\mathfrak{X}$, we see that $M\in\mathfrak{X}$. It means that there is a $\mathfrak{X}$-subgroup $M$  of $G$ such that $H/N=MN/N$. $\square$

\textbf{Lemma 2.2.} \emph{Let $\mathfrak{F}$ be a hereditary saturated formation, $N$ be a normal subgroup of a group $G$, $H$ be a subgroup of $G$     then}

$(1)$  $\mathrm{NI}_\mathfrak{F}(G)N/N\leq \mathrm{NI}_\mathfrak{F}(G/N)$.

$(2)$ $\mathrm{NI}_\mathfrak{F}(G)\cap H\leq  \mathrm{NI}_\mathfrak{F}(H)$.

$(3)$ \emph{Let $N\leq\mathrm{Int}_\mathfrak{F}(G)$ then $N \leq\mathrm{NI}_\mathfrak{F}(G)$ and $\mathrm{NI}_\mathfrak{F}(G)/N=\mathrm{NI}_\mathfrak{F}(G/N)$.}

\emph{Proof.} $(1)$  If $K/N$ is a $\mathfrak{F}$-maximal subgroup of $G/N$ then by lemma 2.1 there exists a $\mathfrak{F}$-maximal subgroup $Q$ of $G$ such that $QN/N=K/N$. If $x\in N_G(Q)$ then $xN\in N_{G/N}(QN/N)=N_{G/N}(K/N)$. Thus  $\mathrm{NI}_\mathfrak{F}(G)N/N\leq \mathrm{NI}_\mathfrak{F}(G/N)$.

$(2)$   If $M$ is a $\mathfrak{F}$-maximal subgroup of $H$ then there exists a $\mathfrak{F}$-maximal subgroup $Q$ of $G$ such that $Q\cap H=M$. So if $x\in \mathrm{NI}_\mathfrak{F}(G)\cap H$  then $M^x=Q^x\cap H^x=Q\cap H=M$. Hence $x\in \mathrm{NI}_\mathfrak{F}(H)$. Thus  $\mathrm{NI}_\mathfrak{F}(G)\cap H\leq  \mathrm{NI}_\mathfrak{F}(H)$.

$(3)$  Let $N\leq\mathrm{Int}_\mathfrak{F}(G)$. It is clear that $N \leq\mathrm{NI}_\mathfrak{F}(G)$. Note that $M$ is a $\mathfrak{F}$-maximal subgroup of $G$ if and only if $M/N$ is a $\mathfrak{F}$-maximal subgroup of $G/N$. Now $N_G(M)/N=N_{G/N}(M/N)$. Thus    $\mathrm{NI}_\mathfrak{F}(G)/N=\mathrm{NI}_\mathfrak{F}(G/N)$. $\square$

Let $\mathfrak{F}$ be a saturated formation. Then in every group exists a $\mathfrak{F}$-projector \cite[p. 292]{s8}.  Recall that a $\mathfrak{F}$-projector  of a group $G$ is a $\mathfrak{F}$-maximal subgroup $H$ of $G$ such that $HN/N$ is a $\mathfrak{F}$-maximal subgroup of $G/N$ for every normal subgroup $N$ of $G$.

Recall that a group $G$ is called semisimple if $G$ is the direct product of simple groups. A chief factor of a group is the example of a semisimple group.

\textbf{Lemma 2.3.} \emph{Let $\mathfrak{F}$ be a hereditary saturated formation and a group $G=HK$ be a product of normal $\mathfrak{F}$-subgroups. If $K$ is semisimple then  $G\in\mathfrak{F}$.}

\emph{Proof.} Assume the contrary. Let a group $G$ be a counterexample of a minimal order. Then $G=HK$ is a product of normal $\mathfrak{F}$-subgroups $H$ and $K$ where $K$ is semisimple. Let $N$ be a normal subgroup of $G$. Then $G/N=(HN/N)(KN/N)$ where $HN/N$ and $KN/N$ are  normal $\mathfrak{F}$-subgroups of $G/N$ and $KN/N$ is semisimple. So $G/N\in\mathfrak{F}$. Since $\mathfrak{F}$ is a saturated formation, we see that $\Phi(G)=1$ and $G$ has an unique minimal normal subgroup that equals $K$. Now $K\leq H$. So $G=H\in\mathfrak{F}$, the contradiction. $\square$


The following lemma is well known.

\textbf{Lemma 2.4.} \emph{Let $\mathfrak{F}$ be a hereditary saturated formation and $H$ be a $\mathfrak{F}$-subgroup of a group $G$. Then $\mathrm{Z}_\mathfrak{F}(G)H\in\mathfrak{F}$. }

Recall that if $\mathfrak{F}$ is a hereditary formation then a subgroup $H$ of a group $G$ is called $\mathfrak{F}$-subnormal if either $H=G$ or there is a chain of subgroups $H=H_0<H_1<\dots< H_n=G$ such that $H_i/\mathrm{Core}_{H_i}(H_{i-1})\in\mathfrak{F}$ for all $i=1, \dots, n.$ We will need the following facts about $\mathfrak{F}$-subnormal subgroups.

\textbf{Lemma 2.5 \cite[p. 236]{s9}.} \emph{Let $\mathfrak{F}$ be a hereditary formation, $N$ be a normal subgroup of a group $G$ and $H$, $K$  be  subgroups of $G$. Then}:

     (1) \emph{If $H$ is $\mathfrak{F}$-subnormal in  $G$ then $HN/N$   is $\mathfrak{F}$-subnormal in $G/N$}.

     (2) \emph{If $H/N$ is $\mathfrak{F}$-subnormal in $G/N$ then $H$   is $\mathfrak{F}$-subnormal in $G$}.

     (3) \emph{If $H$ is $\mathfrak{F}$-subnormal in  $K$  and $K$ is $\mathfrak{F}$-subnormal in  $G$ then $H$   is $\mathfrak{F}$-subnormal in $G$}.

\textbf{Lemma 2.6 \cite[p. 239]{s9}.} \emph{Let $\mathfrak{F}$ be a saturated formation and a group $G=H\mathrm{F}^*(G)$ where $H$ is a $\mathfrak{F}$-subnormal $\mathfrak{F}$-subgroup of $G$. Then $G\in\mathfrak{F}$.  }

\textbf{Lemma 2.7 \cite[p. 390]{s8}.} \emph{Let $\mathfrak{F}$ be a hereditary saturated formation then $[G^\mathfrak{F}, \mathrm{Z}_\mathfrak{F}(G)]=1$ for any group $G$.}


\section{Proves of the main results}

\subsection{Proof of proposition 1}

Let $\mathfrak{F}$ be a hereditary saturated formation and $G$ be a group.   According to lemma 2.2 all  $\mathfrak{F}$-maximal subgroups of $\mathrm{NI}_\mathfrak{F}(G)$ are normal in $\mathrm{NI}_\mathfrak{F}(G)$. Among this $\mathfrak{F}$-maximal subgroups there is a $\mathfrak{F}$-projector $H$. Now  $\mathrm{NI}_\mathfrak{F}(G)/H$ does not contain any $\mathfrak{F}$-subgroup.  Hence $\mathrm{NI}_\mathfrak{F}(G)/H\in\pi(\mathfrak{F})'$.

It is clear that $\mathrm{Int}_\mathfrak{F}(G)\leq\mathrm{O}^{\pi'}(\mathrm{NI}_\mathfrak{F}(G))\in\mathfrak{F}$. Let us show  by induction the the equality $\mathrm{Int}_\mathfrak{F}(G)=\mathrm{O}^{\pi'}(\mathrm{NI}_\mathfrak{F}(G))$ holds. It is clear that it holds for the unit group. Assume that we prove our statement for groups whose order is less then the order of a group $G$.   Let $N$ be a minimal normal subgroup of $G$ such that $N\leq\mathrm{O}^{\pi'}(\mathrm{NI}_\mathfrak{F}(G))$ and   $M$ be a $\mathfrak{F}$-maximal subgroup of $G$. So $M\triangleleft MN$ and $N$ is a normal semisimple $\mathfrak{F}$-subgroup of $MN$. By lemma 2.3 $MN\in\mathfrak{F}$ and hence $MN=M$ for all $\mathfrak{F}$-maximal subgroups $ M$ of $G$. Hence $N\leq \mathrm{Int}_\mathfrak{F}(G)$.

By induction $\mathrm{Int}_\mathfrak{F}(G/N)=\mathrm{O}^{\pi'}(\mathrm{NI}_\mathfrak{F}(G/N))$. According to \cite{h4} $\mathrm{Int}_\mathfrak{F}(G)/N=\mathrm{Int}_\mathfrak{F}(G/N)$. By (3) of lemma 2.2 $\mathrm{NI}_\mathfrak{F}(G)/N=\mathrm{NI}_\mathfrak{F}(G/N)$. From $\pi(N)\subseteq\pi(\mathfrak{F})$ it follows that  $\mathrm{O}^{\pi'}(\mathrm{NI}_\mathfrak{F}(G)/N)=\mathrm{O}^{\pi'}(\mathrm{NI}_\mathfrak{F}(G))/N$. Now $\mathrm{Int}_\mathfrak{F}(G)/N=\mathrm{Int}_\mathfrak{F}(G/N)=\mathrm{O}^{\pi'}(\mathrm{NI}_\mathfrak{F}(G/N))
=\mathrm{O}^{\pi'}(\mathrm{NI}_\mathfrak{F}(G)/N)=\mathrm{O}^{\pi'}(\mathrm{NI}_\mathfrak{F}(G))/N$. Thus $\mathrm{Int}_\mathfrak{F}(G)=\mathrm{O}^{\pi'}(\mathrm{NI}_\mathfrak{F}(G))$. $\square$

\subsection{Proof of theorem A}


The following result directly follows from the proof of the main result of \cite{h4}.

\smallskip

\emph{Let $\mathfrak{H}$ be a hereditary saturated formation and $\pi=\pi(\mathfrak{H})$. Then for every $\pi$-group $G$ the intersection of all $\mathfrak{H}$-maximal subgroups of $G$ is the $\mathfrak{H}$-hypercenter of $G$ if and only if $\mathfrak{H}$ satisfies the boundary condition in the universe of all $\pi$-groups. }

\smallskip

According to \cite[p. 96]{s9} $\mathfrak{F}$ is a hereditary saturated formation. So $\mathfrak{F}$ is a local formation. Let $F$ be the canonical local definition of $\mathfrak{F}$.

$(1)\Rightarrow(2)$ Assume that $\mathfrak{F}_i$  satisfies the boundary condition in the universe of all $\pi_i$-groups for all $i\in I$. Let us show that $\underset{i\in I}\bigcap\mathrm{NI}_{\mathfrak{F}_i}(G)=\mathrm{Z}_\mathfrak{F}(G)$ holds for every group $G$.

Let $G$ be a group and  $D=\underset{i\in I}\bigcap\mathrm{NI}_{\mathfrak{F}_i}(G)$. By proposition 1   $\mathrm{NI}_{\mathfrak{F}_i}(G)$ has the  Hall $\pi_i$-subgroup that belongs to $\mathfrak{F}_i$ and is normal in $G$. Hence $D$ has the normal   Hall $\pi_i$-subgroup that belongs to $\mathfrak{F}_i$ for every $i\in I$. Thus $D\in\mathfrak{F}$.

Let $H/K$ be a chief factor of $G$ below $D$. Then $\pi(H/K)\subseteq\pi_n$ for some $n\in I$.

\smallskip

 (a) $\mathrm{O}^{\pi_n}(G/K)\leq C_{G}(H/K)$.

 By (1) of lemma 2.2 $H/K$ normalizes all $\mathfrak{F}_i$-maximal subgroups of $G/K$.   Hence $H/K$ normalizes all $\mathfrak{F}_i$-projectors of $G/K$ for all $i\in I\setminus\{n\}$.
 Let $F/K$ be a $\mathfrak{F}_i$-projector of $G/K$ for some $i\in I\setminus\{n\}$.
 From $\pi_n\cap\pi_i=\emptyset $  it follows that $F/K\cap H/K=K/K$. Let $hK\in H/K$ and $fK\in F/K$. Then from one hand $[fK, hK]=(fK)^{-1}(fK)^{(hK)}\in F/K$ and from another hand $[fK, hK]=((hK)^{-1})^{(fK)}(hK)\in H/K$. So $[fK, hK]=1$. Hence $[H/K, F/K]=1$.  Thus $H/K$ centralizes all $\mathfrak{F}_i$-projectors of $G/K$ for all $i\in I\setminus\{n\}$.
 Since $\mathfrak{F}_i$ is a hereditary saturated formation for all $i\in I$, we see that $G/C_G(H/K)$ does not contain any $\pi_i$-subgroups for  all $i\in I\setminus\{n\}$. Thus  $\mathrm{O}^{\pi_n}(G/K)\leq C_G(H/K)$.

\smallskip

 (b) \emph{$H/K\leq \mathrm{Z}_{\mathfrak{F}_n}(R/K)$ for every   $\mathfrak{G}_{\pi_n}$-maximal subgroup $R/K$ of $G/K$}.

  Let $R/K$ be a  $\mathfrak{G}_{\pi_n}$-maximal subgroup  of $G/K$. Then $\pi_n((R/K)(H/K))\subseteq\pi_n$. Hence $(R/K)(H/K)=R/K$. So $H/K\subseteq R/K$. By (2) of lemma 2.2 $H/K$ normalizes all $\mathfrak{F}_n$-maximal subgroups of $R/K$. Note that $H/K$ is semisimple $\mathfrak{F}_n$-subgroup. So $(H/K)(F/K)\in\mathfrak{F}_n$ for every $\mathfrak{F}_n$-maximal subgroup $F/K$ of $R/K$ by lemma 2.3. Hence   $(H/K)(F/K)=F/K$ for every $\mathfrak{F}_n$-maximal subgroup $F/K$ of $R/K$. Thus $H/K\leq \mathrm{Int}_{\mathfrak{F}_n}(R/K)$. Since $\mathfrak{F}_n$ satisfies the boundary condition in the universe of all $\pi_n$-groups, $H/K\leq Z_{\mathfrak{F}_n}(R/K)$.

\smallskip

(c)\emph{ Let $R/K$ be a $\mathfrak{G}_{\pi_n}$-maximal subgroup  of $G/K$ such that $(R/K)\mathrm{O}^{\pi_n}(G/K)=G/K$. Then $H/K$ is a chief factor of $R/K.$}

Assume that $N/K$ is a minimal normal subgroup of $R/K$ such that  $K/K\neq N/K<H/K$. From  $\mathrm{O}^{\pi_n}(G/K)\leq C_{G}(H/K)$ it follows that $\mathrm{O}^{\pi_n}(G/K)\leq C_{G}(N/K)$. From  $(R/K)\mathrm{O}^{\pi_n}(G/K)=G/K$ it follows that $N/K$ is normal in $G/K$. Hence $H/K$ is not a chief factor of $G$, a contradiction.

\smallskip

(d)     \emph{$(R/K)^{F(p)}\leq C_G(H/K)$ for all $p\in\pi(H/K)$}.

From $H/K\leq \mathrm{Z}_{\mathfrak{F}_n}(R/K)$ and $\mathfrak{F}_n\subseteq\mathfrak{F}$ it follows that a chief factor $H/K$ of $R/K$ lies in $\mathrm{Z}_\mathfrak{F}(R/K)$. Now $(R/K)/C_{R/K}(H/K)\in F(p)$ for all $p\in\pi(H/K)$. Thus  $(R/K)^{F(p)}\leq C_G(H/K)$ for all $p\in\pi(H/K)$.

\smallskip

(e) \emph{$H/K$ is a $\mathfrak{F}$-central chief factor of $G$}.

From $\mathrm{O}^{\pi_n}(G/K)\leq C_{G}(H/K)$,   $(R/K)^{F(p)}\leq C_G(H/K)$ for all $p\in\pi(H/K)$ and \linebreak$(R/K)\mathrm{O}^{\pi_n}(G/K)=G/K$ it follows that $G/C_G(H/K)\in F(p)$ for all $p\in\pi(H/K)$. Thus $H/K$ is a $\mathfrak{F}$-central chief factor of $G$.

\smallskip

(f) $D\leq \mathrm{Z}_\mathfrak{F}(G)$.

We showed that every chief factor of $G$ below $D$ is $\mathfrak{F}$-central. Hence $D\leq \mathrm{Z}_\mathfrak{F}(G)$.

\smallskip

(g) \emph{$D\geq \mathrm{Z}_\mathfrak{F}(G)$ and hence $D= Z_{\mathfrak{F}}(G)$}.

 Let $H$ be a $\mathfrak{F}_i$-maximal subgroup of $G$ for some $i\in I$. Then $HZ_{\mathfrak{F}}(G)\in \mathfrak{F}$ by lemma 2.4. Since $H$ is  a    $\mathfrak{F}_i$-maximal subgroup of $G$, $H$ is    a  $\mathfrak{F}_i$-maximal subgroup of $HZ_{\mathfrak{F}}(G)$. So $H\triangleleft HZ_{\mathfrak{F}}(G)$. Hence $D\geq Z_{\mathfrak{F}}(G)$. Thus $D= Z_{\mathfrak{F}}(G)$.

\medskip

$(2)\Rightarrow(1)$ Suppose now that $\underset{i\in I}\bigcap\mathrm{NI}_{\mathfrak{F}_i}(G)=\mathrm{Z}_\mathfrak{F}(G)$ holds for every group $G$. Let us show that $\mathfrak{F}_i$  satisfies the boundary condition in the universe of all $\pi_i$-groups for all $i\in I$.

Assume the contrary. Then some $\mathfrak{F}_n$ does not satisfy the boundary condition in the universe of all $\pi_n$-groups. So  there is $\pi_n$-group $G$ such that $\mathrm{Int}_{\mathfrak{F}_n}(G)\neq\mathrm{Z}_{\mathfrak{F}_n}(G)$. Note that $\mathrm{Int}_{\mathfrak{F}_n}(G)=\mathrm{NI}_{\mathfrak{F}_n}(G)$ by proposition 1. Since $G$ is  a $\pi_n$-group,   $\mathrm{NI}_{\mathfrak{F}_n}(G)=\underset{i\in I}\bigcap\mathrm{NI}_{\mathfrak{F}_i}(G)$.
From $\mathfrak{G}_{\pi_n}\cap\mathfrak{F}=\mathfrak{F}_n$ it follows that $\mathrm{Z}_{\mathfrak{F}_n}(G)=\mathrm{Z}_{\mathfrak{F}}(G)$.

Hence $\underset{i\in I}\bigcap\mathrm{NI}_{\mathfrak{F}_i}(G)=\mathrm{NI}_{\mathfrak{F}_n}(G)
=\mathrm{Int}_{\mathfrak{F}_n}(G)\neq\mathrm{Z}_{\mathfrak{F}_n}(G)=\mathrm{Z}_{\mathfrak{F}}(G)$, the contradiction.

\subsection{Proof of proposition 2}

Let $N$ be a normal subgroup of a group $G$ such that $H$ is $\mathfrak{F}$-subnormal in $HN$   for every $H\in\Sigma(G)$. Let $S$ be a $\mathfrak{F}$-subnrmalizer in $G$ of $H\in\Sigma(G)$.  Then $HN/N$ is $\mathfrak{F}$-subnormal in $SN/N$ by (1) of lemma 2.5. So  $HN$ is $\mathfrak{F}$-subnormal in $SN$ by (2) of lemma 2.5. Hence $H$ is $\mathfrak{F}$-subnormal in $SN$ by (3) of lemma 2.5. Thus $SN=N$. It means that $N\leq \mathrm{SI}_\Sigma^\mathfrak{F}(G).$ So every normal subgroup of $G$ that $\mathfrak{F}$-subnormalize all subgroups from $\Sigma(G)$ lies in   $\mathrm{SI}_\Sigma^\mathfrak{F}(G).$

From the other hand $H\mathrm{SI}_\Sigma^\mathfrak{F}(G)$ belongs to every $\mathfrak{F}$-subnormalizer of $H$ in $G$   for every $H\in\Sigma(G)$. Hence $H$ is $\mathfrak{F}$-subnormal in $H\mathrm{SI}_\Sigma^\mathfrak{F}(G)$ for every $H\in\Sigma(G)$.

\subsection{Proof of theorem B}

$(1)\Rightarrow(2)$ Assume that  there exists a partition $\sigma=\{\pi_i | i\in I\}$ of $\pi(\mathfrak{F})$ into disjoint subsets such  that $\mathfrak{F}=\underset{i\in I}\times\mathfrak{G}_{\pi_i}$. Let us show that the intersection of all $\mathfrak{F}$-subnormalizers of all cyclic primary subgroups of a group $G$ is the $\mathfrak{F}$-hypercenter of $G$ for every group $G$.

Note that $\mathfrak{F}$ is local formation with the canonical local definition $F$ where $F(p)=\mathfrak{G}_{\pi_i}$ for $p\in\pi_i$ for all $i\in I$.

Let $D$ be the intersection of all $\mathfrak{F}$-subnormalizers of all cyclic primary subgroups of a group $G$ and $H/K$ be a chief factor of $G$ below $D$.

\smallskip

(a) \emph{$H/K$ lies in the intersection of all $\mathfrak{F}$-subnormalizers of all cyclic primary subgroups of a group $G/K$}.

Let $CK/K$ be a cyclic primary subgroup of $G/K$.  According to lemma 2.1 we may assume that $C$ is    a cyclic primary subgroup of $G$. Now $C$ is $\mathfrak{F}$-subnormal in $HC$ by proposition 2. So $CK/K$ is      $\mathfrak{F}$-subnormal in $HC/K$ by (1) of lemma 2.5. Hence $H/K$ lies in the intersection of $\mathfrak{F}$-subnormalizers of all cyclic primary subgroups of $G/K$.

\smallskip

(b) $H/K\in\mathfrak{F}$.

Now $K/K$ is a $\mathfrak{F}$-subnormal $\mathfrak{F}$-subgroup of a quasinilpotent group $H/K$. By lemma 2.6 $H/K\in\mathfrak{F}$. Hence $\pi(H/K)\subseteq\pi_n$ for some $n\in I$.

\smallskip

(c) \emph{$C/K\leq C_G(H/K)$ for every cyclic primary $\pi(\mathfrak{F})'$-subgroup of $G/K$}.

Let $C/K$ be a cyclic primary $\pi(\mathfrak{F})'$-subgroup of $G/K$. Since $C/K$ is a $\mathfrak{F}$-subnormal $\pi(\mathfrak{F})'$-subgroup of $HC/K$, $(C/K)^\mathfrak{F}=(C/K)$ is subnormal in $HC/K$ by (1) of lemma 6.1.9 \cite[p. 237]{s9}. Since $C/K$ is a subnormal Sylow subgroup of $HC/K$, we see that $C/K\triangleleft HC/K$. Now it is easy to see that $C/K\leq C_G(H/K)$.

\smallskip

(d) \emph{$C/K\leq C_G(H/K)$ for every  cyclic primary $\pi(\mathfrak{F})\cap(\pi_n')$-subgroup of $G/K$}.

Let $C/K$ be a cyclic primary $\pi(\mathfrak{F})\cap(\pi_n')$-subgroup of $G/K$. Since $C/K$ is $\mathfrak{F}$-subnormal in $HC/K$, $HC/K\in\mathfrak{F}$  by lemma 2.6. So $C/K\leq C_G(H/K)$.

\smallskip

(e) \emph{$H/K$  is a $\mathfrak{F}$-central chief factor of $G$ and $D\leq \mathrm{Z}_\mathfrak{F}(G)$}.

From (c) and (d) it follows that $\mathrm{O}^{\pi_n}(G)\leq C_G(H/K)$. Hence $G/C_G(H/K)\in \mathfrak{G}_{\pi_n}=F(p)$ for all $p\in\pi(H/K)$. So $H/K$  is a $\mathfrak{F}$-central chief factor of $G$. It means that $D\leq \mathrm{Z}_\mathfrak{F}(G)$.

\smallskip

(f) \emph{$\mathrm{Z}_\mathfrak{F}(G)\leq D$ and hence $D=\mathrm{Z}_\mathfrak{F}(G)$}.

Let $C$ be a cyclic $p$-subgroup of a group $G$. If $p\in\pi(\mathfrak{F})$ then $C\mathrm{Z}_\mathfrak{F}(G)\in\mathfrak{F}$ by lemma 2.4. Hence $C$ is $\mathfrak{F}$-subnormal in $C\mathrm{Z}_\mathfrak{F}(G)$.  If $p\not\in\pi(\mathfrak{F})$ then $C\leq G^\mathfrak{F}$. By lemma 2.7 $C\leq C_G(\mathrm{Z}_\mathfrak{F}(G))$. Hence $(C\mathrm{Z}_\mathfrak{F}(G))^\mathfrak{F}=C$. So $C$ is $\mathfrak{F}$-subnormal in $C\mathrm{Z}_\mathfrak{F}(G)$. Hence $\mathrm{Z}_\mathfrak{F}(G)\leq D$. Thus $D=\mathrm{Z}_\mathfrak{F}(G)$.

$(2)\Rightarrow(3)$ Let $P$ be a Sylow $p$-subgroup of $G$. If  $p\in\pi(\mathfrak{F})$ then $P\in\mathfrak{F}$ and hence $P\mathrm{Z}_\mathfrak{F}(G)\in\mathfrak{F}$ by lemma 2.4. So $P$ is $\mathfrak{F}$-subnormal in $P\mathrm{Z}_\mathfrak{F}(G)$.

 If  $p\not\in\pi(\mathfrak{F})$ then $P\leq G^\mathfrak{F}$. By lemma 2.7 $[G^\mathfrak{F}, \mathrm{Z}_\mathfrak{F}(G)]=1$. So $P\mathrm{Z}_\mathfrak{F}(G)=P\times\mathrm{Z}_\mathfrak{F}(G)$.     Hence $P$ is $\mathfrak{F}$-subnormal in $P\mathrm{Z}_\mathfrak{F}(G)$.

 Thus $\mathrm{Z}_\mathfrak{F}(G)$ lies in the intersection $D$ of all $\mathfrak{F}$-subnormalizers of all Sylow subgroups of $G$. Since the unit subgroup is $\mathfrak{F}$-subnormal in $D$, we see that $\pi(D)\subseteq\pi(\mathfrak{F}).$

 Now let $C$ be a cyclic primary $p$-subgroup of $G$. Then there is a Sylow $p$-subgroup $P$ of $G$ such that $C\leq P$.     If   $p\in\pi(\mathfrak{F})$ then $C$ is $\mathfrak{F}$-subnormal in $P$ and $P$ is $\mathfrak{F}$-subnormal in $PD$. Hence $C$ is $\mathfrak{F}$-subnormal in $PD$ and also in $CD$ by lemma 2.5.

   If   $p\not\in\pi(\mathfrak{F})$ then $C$ is subnormal in $P$ and $P$ is normal in $PD$. Hence $C$ is subnormal in $PD$ and also in $CD$. So $C$ is the normal Sylow subgroup of $CD$. By our assumption the unit group is a Sylow subgroup. Hence $1$ is $\mathfrak{F}$-subnormal in $D$. Now  $C/C$ is $\mathfrak{F}$-subnormal in $CD/C$. Hence $C$ is $\mathfrak{F}$-subnormal in $CD$.

   Thus $D$ lies in the intersection of all $\mathfrak{F}$-subnormalizers of all cyclic primary subgroups of $G$. Hence $\mathrm{Z}_\mathfrak{F}(G)\leq D\leq \mathrm{Z}_\mathfrak{F}(G)$. Thus $D=\mathrm{Z}_\mathfrak{F}(G)$.

Consider the following statement:

(4) \emph{$\mathfrak{F}$ has the canonical local definition $F$ such that for every prime $p$, $F(p)$ contains every group $G$ whose all Sylow subgroups belong to $F(p)$}.

$(3)\Rightarrow(4)$ Let the intersection of all $\mathfrak{F}$-subnormalizers of all Sylow subgroups of   $G$ be the $\mathfrak{F}$-hypercenter of $G$ for every group $G$. Assume that there exist a prime $p$ and  groups $G$ such that $G\not\in F(p)$ but for every Sylow subgroup $P$ of $G$, $P\in F(p)$. Let us chose the minimal order group $G$ from such groups.

It is clear that $\mathrm{O}_p(G)=1$ and $G$ has an unique minimal normal subgroup. Then by lemma 2.6 from \cite{h4} there exists a faithful irreducible $G$-module $N$ over the field $\mathrm{F}_p$. Let $H$ be the semidirect product of $N$ and $G$. Note that $NP\in\mathfrak{F}$ for every Sylow subgroup $P$ of $H$. Hence $N$ lies in the intersection of all $\mathfrak{F}$-subnormalizers of Sylow subgroups of $H$ by proposition 2. But $H/C_H(N)\not\in F(p)$. So $N\not\leq \mathrm{Z}_\mathfrak{F}(H)$,   the contradiction.

$(4)\Rightarrow (1)$. Assume that $Z_q\in F(p)$ for primes $p\neq q$. Suppose that $F(q)\cap \mathfrak{N}_p\neq\mathfrak{N}_p$.  Let $P$ be the minimal order $p$-group from $\mathfrak{N}_p\setminus (F(q)\cap \mathfrak{N}_p)$. Then $P$  has an unique minimal normal subgroup and $P\in F(p)$. There exists a faithful irreducible $P$-module $Q$ over the field $\mathrm{F}_q$. Note that $Q\in F(p)$. Hence the semidirect product $G=Q\leftthreetimes P\in F(p)\subseteq \mathfrak{F}$. Now $G/\mathrm{O}_{q', q}(G)= G/Q\simeq P\in F(q)$, a contradiction.

So from $Z_q\in F(p)$ it follows that $F(q)\cap \mathfrak{N}_p=\mathfrak{N}_p$ and hence $F(p)\cap \mathfrak{N}_q=\mathfrak{N}_q$. So $\mathfrak{N}_{\pi(F(p))}\subseteq F(p)$. Let a group $G$ be a $s$-critical for $F(p)$. Since $F(p)$ is hereditary, we see that $G$ is $r$-group for some prime $r$. Now $r\not\in\pi(F(p))$. Hence $G\simeq Z_r$. It means that $F(p)=\mathfrak{G}_{\pi(F(p))}$ for all $p\in\pi(\mathfrak{F})$.

Assume   now that for three different primes $p, q$ and $r$ we have that $\{p, q\}\subseteq\pi(F(r))$. Let us show that $q\in\pi(F(p))$. By theorem 10.3B \cite{s8} there exists a faithful irreducible $Z_q$-module $P$ over the field $\mathrm{F}_p$. Let $G$ be the semidirect product $P$ and $Z_q$. Then $T\in F(r)\subseteq \mathfrak{F}$.  Thus $G/\mathrm{O}_{p', p}(G)= G/P\simeq Z_q\in F(p)$.

It means that there exists a partition $\sigma=\{\pi_i | i\in I\}$ of $\pi(\mathfrak{F})$ into disjoint subsets such that $F(p)=\mathfrak{G}_{\pi_i}$ for all $p\in\pi_i$ and for all $i\in I$. Now    $\mathfrak{F}=\underset{i\in I}\times\mathfrak{G}_{\pi_i}$.



\end{document}